\documentclass{amsart}

\usepackage[all]{xy}
\usepackage{amsmath}\usepackage{amssymb}

\newcommand{\DG}[1]{\cite[{#1}]{mdgg:crahcrrz}}
\newcommand{\DL}[1]{\cite[{#1}]{mdtal:ccomk}}
\newcommand{\DandG}{\cite{mdgg:crahcrrz}}
\newcommand{\DandLii}{\cite{mdtal:umctkg}}
\newcommand{\DandL}{\cite{mdtal:ccomk}}

\newcommand{\KK}{K\!K}
\newcommand{\Ks}{K_\star}
\newcommand{\Marius}{D{\u{a}}d{\u{a}}rlat}
\newcommand{\NN}{{\Bbb{N}}}
\newcommand{\QQ}{{\Bbb{Q}}}

\newcommand{\ZZp}[1]{\ZZ/{#1}}
\newcommand{\ZZ}{{\Bbb{Z}}}
\newcommand{\bmapt}{\tilde{\beta}}
\newcommand{\bmap}{\beta}

\newcommand{\cstar}{\mbox{$C^*$}}

\newcommand{\gritem}[1]{\noindent{\sl{#1}$^\circ$:}}
\newcommand{\id}{\operatorname{id}}
\newcommand{\indlim}{{\displaystyle\lim_{\longrightarrow}}}
\newcommand{\infer}{\mbox{$\Longrightarrow$}}
\newcommand{\inv}{^{-1}}

\newcommand{\kmapk}[2]{\kappa_{{#1},{#2}}}

\newcommand{\lmapk}[2]{\lambda_{{#1},{#2}}}

\newcommand{\mydiv}{\uparrow}
\newcommand{\mythesis}[1]{\cite[\ref{#1}]{se:iaa}}

\newcommand{\pKsubi}[2]{K_i({#2};\ZZp{#1})}

\newcommand{\pKu}[2]{K_0({#2};\ZZ\oplus\ZZp{#1})}
\newcommand{\pKzm}[1]{K_0({#1};\QQ/\ZZ)}
\newcommand{\pK}[2]{K_0({#2};\ZZp{#1})}
\newcommand{\po}{^+}
\newcommand{\restr}[1]{\left|{\raisebox{-.7ex}[.7ex][.7ex]{${#1}$}}\right.}
\newcommand{\rmapt}{\tilde{\rho}}
\newcommand{\rmap}{\rho}
\newcommand{\scale}{\Sigma}

\newcommand{\tor}{\operatorname{tor}}

\newcommand{\PK}[1]{\overline{\mathbf K}(#1)}
\newcommand{\pKz}[1]{K_0(#1;\QQ/\ZZ)}

\newcounter{rocount}

\newenvironment{gralist}{\begin{list}{\arabic{rocount}$^\circ$:}{\usecounter{rocount}}}{\end{list}}

\newtheorem{theorem}{Theorem}[section]
\newtheorem{proposition}[theorem]{Proposition}
\newtheorem{lemma}[theorem]{Lemma}
\newtheorem{definition}[theorem]{Definition}
\newtheorem{remark}[theorem]{Remark}

\newenvironment{pf}{{\bf\sl Proof:}}{\hfill$\blacksquare$}

\begin{document}

\title
{K\"unneth Splittings and Classification of \cstar-Algebras with
Finitely many Ideals}
\author{S\o ren Eilers}

\begin{abstract}
The class of $AD$ algebras of real rank zero is classified by an exact
sequence of $K$-groups with coefficients, equipped with certain
order structures. Such a sequence is always {\em split}, and one may
ask why, then, the middle group is relevant for classification. The
answer is that the splitting map can not always be chosen to respect
the order structures involved.

This may be rephrased in terms of the ideals of the \cstar-algebras in
question. We prove that when
the \cstar-algebra has  only finitely many ideals, a splitting map
respecting these always exists. Hence $AD$ algebras of real rank zero
with finitely many ideals are classified by (classical) ordered
$K$-theory. We also indicate how the methods generalize to the full
class of $ASH$ algebras with slow dimension growth and real rank zero. 
\end{abstract}
\maketitle
\section{Introduction}
In the ongoing quest to classify \cstar-algebras of real rank zero by
algebraic invariants much progress has been made on classes of
\cstar-algebras defined from a narrow class of {\em building blocks}\/
as the smallest collection, containing these, which is closed under
taking matrices, finite direct sums, and countable inductive limits.
Considering the single building block $C(S^1)$ one gets the $AT$
algebras, whereas the more extensive $AH$ class is obtained by
allowing every commutative \cstar-algebra $C(Y)$ for $Y$ 
a finite $CW$-complex. The $AD$ algebras are derived from the building
blocks $C(S^1)$ and the dimension drop intervals $\Bbb I_n^{\sim}$, and the
$ASH$ algebras are defined allowing both all finite $CW$ complexes and
the dimension drop intervals. In the cases of classes where the
dimension of the topological spaces is allowed to vary ($AH$, $ASH$)
one must often impose the {\it slow dimension growth} condition
introduced in \cite{bbmdmr:rrilca}.

In recent years, examples have appeared to
demonstrate that ordered, graded $K$-theory is {\em not}\/ a complete
invariant for any class of stable rank one, real rank zero
\cstar-algebras with slow dimension growth that extends much beyond the $AT$ class, unless one
imposes restrictions on the ideal structure of the algebras in
question.  The first such example -- a pair of nonisomorphic
$AH$ algebras with slow dimension growth having isomorphic ordered, graded $K$-theory --  was found by Gong
in \cite{gg:ccrrzuet}. Subsequently, Elliott-Gong-Su in 
\cite[2.19]{egs:ccrrzrlsdt} provided similar examples involving algebras which
were simultaneously $AH$ and $AD$, and 
\Marius-Loring (\DandL) found a pair of $AD$ algebras of real rank
zero of which one was $AH$ with slow dimension growth and the other
was not, (cf. \DG{10.21}). 

These examples extinguished the hope, lit by results of Zhang, that
many stable rank one, 
real rank zero \cstar-algebras would be classified by this invariant,
regardless of 
ideals. Indeed, it was proven in \cite{sz:rdpisma}
that the ideal lattice of a \cstar-algebra $A$ in this class is
reflected by the ideal lattice of $K_0(A)$, so as an order isomorphism of
$K_0$-groups hence preserves ideals, it was believed that
by extending this order to the graded group $K_0(A)\oplus K_1(A)$, one
could achieve that an order isomorphism $(\varphi_0,\varphi_1)$ was induced
by a $*$-isomorphism respecting the pairing of ideals given by
$\varphi_0$. This turned out not to be the case in general.

Several augmented invariants have been introduced recently to remedy
this situation (\cite{se:ciabtk}, \DandLii, \DandG), all
based on the 
ingenious definition of an order structure on $K$-theory with
coefficients introduced in \DandL. That order structure was defined
by $\KK$-theoretical means, but it has subsequently proven useful to
work instead with the classical order structures on certain tensor
products as suggested already in \cite[5.6]{bb:stc}.

These invariants are complete for
large classes of stably finite, real rank zero \cstar-algebras which
may have arbitrary ideal lattices. It is the purpose of this note to
show how, in the case that the ideal lattices of the \cstar-algebras to
be classified are {\em finite}, one may induce from an order isomorphism
of the  graded $K$-groups an isomorphism of the larger
invariants, hence proving completeness of the classical invariant in
this case.

More specifically, we shall prove:

\begin{theorem}\label{finite}
The invariant
\[
[\Ks(-),\Ks(-)\po,\scale(-)]
\]
is complete for the class of $AD$ algebras with real rank zero and
finitely many ideals.
\end{theorem}

The first result on this form was given by Gong in \cite{gg:ccrrzuet},
and the idea of the proof given there surely may be employed in our
setting also. Conversely, with a little more work our line of proof
carries over to the invariants considered in \DandG, and hence covers
the full  
case of $ASH$ algebras of real rank zero and slow dimension growth. 
Rather than giving the most general classification result, it is our prime
objective to display a fundamental phenomenon that accounts 
for the fact that 
the generalized $K$-groups are irrelevant in the finite ideal lattice
case. 
This comes out most clearly when one works with one of the subclasses
of the $ASH$ algebras which are classified by a finite number of
$K$-groups. We have chosen to work with the $AD$ algebras, which are
classified by a single short exact sequence of (ordered) $K$-groups,
and are rewarded for making this choice by a number of shortcuts along
the way. To substantiate our claims about the general case we briefly
indicate at the end of the paper what must be done here.

\section{Invariants for $AD$ algebras}

The invariant $\mathbf K(-;n)$ defined in \cite{se:ciabtk} consists of
the ordered groups $\Ks(A)$ and $\pKu{n}{A}$ (see \DandL) along with the maps at
the center of 
the exact sequence
$$
\xymatrix{
{K_0(A)}\ar[r]^-{\times n}&{K_0(A)}\ar[r]^-{\rmap_n}&{\pK{n}{A}}\ar[r]^-{\bmap_n}&
{K_1(A)}\ar[r]{\times n}&{K_1(A)}.}
$$
There are also natural maps
$$
\kmapk{m}{n}:\pK{n}{A}\rightarrow \pK{m}{A},
$$
satisfying
\begin{eqnarray}\label{bkcoh}
\bmap_m\kmapk{m}{n}&=&\tfrac{n}{(n,m)}\bmap_n\\
\label{rkcoh}\kmapk{m}{n}\rmap_n&=&\tfrac{m}{(n,m)}\rmap_m\\
\label{kkcoh}\kmapk{k}{m}\kmapk{m}{n}&=&\tfrac{m(k,n)}{(k,m)(m,n)}\kmapk{k}{n}.
\end{eqnarray}
After taking inductive limits over $\NN$ ordered by
$$
n\leq m\Longleftrightarrow n\mydiv m
$$
we get
$$
\xymatrix{
{K_0(A)}\ar[r]^-{\id\otimes 1}&{K_0(A)\otimes\QQ}\ar[r]^-{\rmap}&{\pKz{A}}\ar[r]^-{\bmap}&
{K_1(A)}\ar[r]^-{\id\otimes 1}&{K_1(A)\otimes\QQ}.}
$$ 
where
$$
\pKz{A}=\indlim{(\pK{n}{A},\kmapk{kn}{n})}.
$$
The orders on $\pKu{n}{A}$ induce one on $K_0(A)\otimes\QQ\oplus\pKz
{A}$, and we get an invariant $\PK{A}$. We then have
(\cite[5.3]{se:ciabtk}, \mythesis{mainGENclas}):

\begin{theorem}\label{clas}
$\PK{-}$ is a complete
invariant for the class of $AD$ algebras of real rank zero.
$\mathbf{K}(-;n)$ is a complete invariant for the subclass where $n$
annihilates the torsion of $K_1(-)$.
\end{theorem}

As these invariants include the maps $\rmap,\bmap$, isomorphisms of
the invariants are triples of positive isomorphisms which are also 
complex homomorphisms, i.e., commute with the maps in the complex.
In the case of $\PK{-}$, we shall require that the map defined on
$K_0(-)\otimes \QQ$ is the one induced by the map $\varphi_0$ on
$K_0(-)$, i.e., of the form $\varphi_0\otimes\id$. As a consequence of
torsion freeness, this is equivalent to requiring the maps to commute
with $\id\otimes -$.

\section{Ideal K\"unneth splittings}
Unsplicing the complexes discussed above, we get {\em K\"unneth
sequences}
\[
\xymatrix{
{0}\ar[r]&{K_0(A)\otimes\ZZp{n}}\ar[r]^-{\rmapt_n}&{\pK{n}
{A}}\ar[r]^-{\bmapt_n}& {K_1(A)[n]}\ar[r]& 0\\
{0}\ar[r]&{K_0(A)\otimes\QQ/\ZZ}\ar[r]^-{\rmapt}&
{\pKz
{A}}\ar[r]^-{\bmapt}& {\tor(K_1(A))}\ar[r]& 0. }
\]
It is clear by injectivity of
$K_0(A)\otimes \QQ/\ZZ$ that the second sequence splits,
and by results of B\"odigheimer (\cite{cfb:skskI}), so does the first.
In the second case, we shall refer to a splitting map
\[
\sigma:\tor(K_1(A))\rightarrow \pKz{A}
\]
as a {\em K\"unneth splitting}. See Remark \ref{ASH} for what should
be understood as a K\"unneth splitting in the first case.

Let $A$ be an $AD$ algebra of real rank zero and let $I$ be an ideal
of $A$. It is well known (\cite{mdtal:ecrrzc}) that both $I$ and $A/I$
are $AD$ algebras of real rank zero and that the six term exact
sequence in $K$-theory then becomes \[ \xymatrix{
0\longrightarrow{\Ks(I)}\longrightarrow{\Ks(A)}\longrightarrow{\Ks(A/I)}\longrightarrow 0}  \] Furthermore,
both sequences are {\em pure exact} by \cite[8,11]{lgbmd:ecq}, in the
case of $K_0$, simply because all groups are torsion free. We get that
the vertical maps in
\begin{flushleft}
\begin{picture}(0,90)
\put(00,70){$0$}
\put(10,70){$\vector(1,0){30}$}
\put(45,70){${K_0(I)\otimes\QQ/\ZZ}$}
\put(110,70){$\vector(1,0){30}$}
\put(120,75){$\tilde{\rmap}$}
\put(70,60){$\vector(0,-1){30}$}
\put(145,70){${\pKz {I}}$}
\put(200,70){$\vector(1,0){30}$}
\put(210,75){${\tilde{\bmap}}$}
\put(170,60){$\vector(0,-1){30}$}
\put(240,70){${\tor(K_1(I))}$}
\put(300,70){$\vector(1,0){30}$}
\put(340,70){$0$}
\put(00,10){$0$}
\put(10,10){$\vector(1,0){30}$}
\put(45,10){${K_0(A)\otimes\QQ/\ZZ}$}
\put(110,10){$\vector(1,0){30}$}
\put(120,15){${\tilde{\rmap}}$} 
\put(145,10){${\pKz{A}}$}
\put(200,10){$\vector(1,0){30}$}
\put(260,60){$\vector(0,-1){30}$}
\put(210,15){${\tilde{\bmap}}$}
\put(240,10){${\tor(K_1(A))}$}
\put(300,10){$\vector(1,0){30}$}
\put(340,10){$0$}
\end{picture} 
\end{flushleft}
are all
embeddings. This is a consequence of torsion freeness of $K_0(A/I)$ to
the left, and in the middle then follows from exactness. We will
identify the $K$-groups derived from $I$ with their images in the
$K$-groups derived from $A$.

We say that a complex homomorphism 
$\Phi=(\varphi_0\otimes\id,\varphi,\varphi_1)$ from $\PK{A}$ to $\PK{B})$ is {\em
ideal-preserving} when $\varphi_0 
K_0(I)\subseteq K_0(J)$ implies $\Phi\PK{I}\subseteq
\PK{J}$, i.e.\
\[
\varphi\pKz{I}\subseteq \pKz{J}\qquad
\varphi_1K_1(I)\subseteq K_1(J).
\]
An {\em ideal isomorphism}\/ is an isomorphism which is
ideal-preserving and has ideal-preserving inverse. This concept
replaces a similar notion for $\KK$-theory or $E$-theory introduced in
\cite{gg:ccrrzuet}. 

The following
result, explaining the relevance of the order structures when working
with ideals, is essentially contained in \DG{4.12,9.2}. A detailed
alternative proof will appear in \cite{mdse:ccpikc}.

\begin{proposition}\label{kuidealprespos}
  Let $A,B$ be real rank zero $AD$ algebras and assume that
  $\Phi=(\varphi_0\otimes\id,\varphi,\varphi_1)$ is an isomorphism of the
  complexes $\PK{A}$ and $\PK{B}$. Then the following conditions are
  equivalent:
  \begin{itemize}
  \item[(i)] $\Phi$ is an order isomorphism
  \item[(ii)] $\varphi_0$ is an order isomorphism and $\Phi$ is an ideal
    isomorphism.
  \end{itemize}
\end{proposition}

Because $\rmap,\bmap$ are natural maps, they always preserve ideals;
for instance, $\rmap(K_0(I))$ $\subseteq \pK{n}{I}$. K\"unneth
splittings are not natural, and so we need the following definition.

\begin{definition}\label{idspli}
  A \cstar-algebra $A$ is {\em ideally split}\/ if there is a K\"unneth
  splitting $\sigma$ preserving ideals, i.e., \[ \sigma
  \tor(K_1(I))\subseteq \pKz{I} \] for all ideals $I$ of $A$.
\end{definition}

\begin{remark}\label{exx}\mbox{}
  \begin{gralist}
  \item A simple $AD$ algebra is ideally split.
  \item An $AD$ algebra with torsion free $K_1$ is ideally split.
  \item  An $AD$ algebra with divisible $K_0$ is ideally split. For
    then $\bmapt$ is invertible, and $\bmapt\inv$ provides a natural,
    hence ideal-preserving, splitting map.
  \item An $AD$ algebra for which every proper ideal has torsion
    free $K_1$ is ideally split. For as the map from $K_1(I)$ to
    $K_1(A)$ is an embedding, the image of $K_1(I)$ misses $\tor
    K_1(A)$.
  \item   Not all $AD$ algebras of real rank zero are ideally split. 
    Consider for instance the algebra $D_p$ in \cite[3.3]{mdtal:ccomk}. 
    We have $K_1(D_p)=\ZZp{p}$, and \[
    \pK{p}{D_p}=\left\{(\overline{a},\overline{b},\overline{c}_i)
    \in\ZZp{p}\oplus\ZZp{p}\oplus \prod_{\ZZ}{\ZZp{p}} \left|
    \begin{array}{l}
      \overline{c}_i=\overline{a}\text{ as }i\rightarrow
      \infty\\ \overline{c}_i=\overline{b}\text{ as }i\rightarrow
      -\infty
    \end{array}\right\}\right.
    \] There are ideals $I_n$ with $K_1(I_n)=\ZZp{p}$ and
    $\pK{p}{I_n}=\{(a,b,c_i)\mid c_i=0, |i|\leq n\}$, so an ideal
    splitting must map into $\ZZp{p}\oplus\ZZp{p}\oplus 0$. As this
    set intersects trivially with $\pK{p}{D_p}$, there is no ideal
    splitting.  (In fact, the algebra $C_p$ considered in \DL{3.3}
    {\em is}\/ ideally split).
  \end{gralist}
\end{remark}

The reader has probably already guessed why we are interested in
ideally split algebras. Here is the accurate statement.

\begin{proposition}\label{idealsplit}
  If $A$ and $B$ are ideally split $AD$ algebras of real rank zero, and
  \[ [\Ks(A),\Ks(A)\po,\scale(A)]\simeq [\Ks(B),\Ks(B)\po,\scale(B)] \]
  then $A\simeq B$.
\end{proposition}
\begin{pf}
  Let $A,B$ be ideally split $AD$ algebras of real rank zero, and let
  $(\varphi_0,\varphi_1)$ be an order isomorphism of $\Ks(A)$ with $\Ks(B)$.
  There is a diagram
\begin{flushleft}
\begin{picture}(0,90)
\put(00,70){$0$}
\put(10,70){$\vector(1,0){30}$}
\put(45,70){${K_0(A)\otimes\QQ/\ZZ}$}
\put(110,70){$\vector(1,0){30}$}
\put(120,75){${\tilde{\rho}}$}
\put(70,60){$\vector(0,-1){30}$}
\put(75,50){${\varphi_0}$}
\put(145,70){${K_0(A,\QQ/\ZZ)}$}
\put(200,70){$\vector(1,0){30}$}
\put(200,63){$\prec\ldots\ldots$}
\put(210,75){${\tilde{\bmap}}$}
\put(210,55){${\sigma}$}
\put(170,60){$\vector(0,-1){30}$}
\put(175,50){$\varphi$}
\put(240,70){${\tor(K_1(A))}$}
\put(260,60){$\vector(0,-1){30}$}
\put(265,50){$\varphi_1$}
\put(300,70){$\vector(1,0){30}$}
\put(340,70){$0$}
\put(00,10){$0$}
\put(10,10){$\vector(1,0){30}$}
\put(45,10){${K_0(B)\otimes\QQ/\ZZ}$}
\put(110,10){$\vector(1,0){30}$}
\put(120,15){${\tilde{\rho}}$}
\put(145,10){${K_0(B,\QQ/\ZZ)}$}
\put(200,10){$\vector(1,0){30}$}
\put(200,05){$\prec\ldots\ldots$}
\put(210,15){${\tilde{\bmap}}$}
\put(210,00){$\tau$}
\put(240,10){${\tor(K_1(B))}$}
\put(300,10){$\vector(1,0){30}$}
\put(340,10){$0$} 
\end{picture}
\end{flushleft}
 in which
  $\sigma$ and $\tau$ are ideal splittings  and $\varphi$ is induced by \[
  \varphi(\tilde{\rmap}(x)+\sigma(y))=\tilde{\rmap}(\varphi_0(x))+\tau(\varphi_1(y)).
  \] As when proving Proposition  \ref{kuidealprespos}, one gets that since $(\varphi_0,\varphi_1)$ is an order
  isomorphism, we have \[ \varphi_0K_0(I)\subseteq K_0(J)\infer
  \varphi_1K_1(I)\subseteq K_1(J) \] (and similarly for
  $(\varphi_0\inv,\varphi_1\inv)$). We hence only need to prove that $\varphi$
(and $\varphi\inv$)
  preserves ideals, and by its definition, this follows by the
  ideal-preserving properties of $\varphi_1,\sigma,$ and $\tau$. The
  triple 
  $(\varphi_0,\varphi,\varphi_1)$ will be an order isomorphism by
  Proposition \ref{kuidealprespos}, and  we may apply Theorem \ref{clas}.
\end{pf}

\begin{remark}\label{someclas}
  Combining Proposition  \ref{idealsplit} with Remark \ref{exx} $1^\circ$--$2^\circ$ we
  regain the well-known classification results, essentially contained
in   \cite{gae:ccrrz}, that $AD$ algebras which are either simple or
have torsion free $K_1$ are classified by their ordered, graded
$K$-groups. See \DL{4.2} and \cite[5.4]{se:ciabtk}. 
\end{remark}

\section{Building ideal splittings}

In this section, we shall prove

\begin{proposition}\label{existidspli}
Any $AD$ algebra of real rank zero with finitely many ideals is
ideally split.
\end{proposition}

Combining this with  Proposition  \ref{idealsplit}, we get Theorem
\ref{finite}. Note how the arguments are predominantly algebraic.
We only use \cstar-algebra results to get purity of certain exact
sequences as mentioned above, and to conclude that a certain lattice of
subgroups is distributive, owing to the fact that the lattice of
ideals of a \cstar-algebra has this property. Such a result can also,
as in \cite{gg:ccrrzuet}, be obtained by appealing to the inductive
limit structure of the \cstar-algebras in question.

\begin{lemma}\label{maximalcase}
Let $I$ be an ideal of a real rank zero $AD$ algebra $A$. Suppose a
splitting map
$$\xymatrix{
0\ar[r]&
{K_0(I)\otimes\QQ/\ZZ}\ar[r]^-{\tilde{\rmap}_I}&
{\pKzm{I}}\ar[r]^-{\tilde{\bmap}_I}&
{\tor(K_1(I))}\ar[r]\ar@{..>}@/_5mm/[l]_{\tau} &0}
$$
is given. Then there is a splitting map
$$\xymatrix{
0\ar[r]&
{K_0(A)\otimes\QQ/\ZZ}\ar[r]^-{\tilde{\rmap}_A}&
{\pKzm{A}}\ar[r]^-{\tilde{\bmap}_A}&
{\tor(K_1(A))}\ar[r]\ar@{..>}@/_5mm/[l]_{\sigma} &0}
$$
extending $\tau$.
\end{lemma}
\begin{pf}
Choose, by Zorn's lemma, a subgroup $D$ of $\pKz{A}$ maximal with
respect to the properties
$$
D\cap {\operatorname{im}}\rmapt_A=0\qquad {\operatorname{im}}\tau\subseteq D.
$$
By maximality, $({\operatorname{im}}\rmapt_A+D)/D$ is an essential subgroup in $\pKz
{A}/D$, and we have
$$
0\longrightarrow {\operatorname{im}}\rmapt_A\longrightarrow \pKz{A}/D\longrightarrow
\dfrac{\pKz{A}/D}{({\operatorname{im}}\rmapt_A+D)/D}\longrightarrow 0
$$
where injectivity to the left is a consequence of
$D\cap{\operatorname{im}}\rmapt_A=0$. This sequence splits by divisibility of
${\operatorname{im}}\rmapt_A\simeq K_0(A)\otimes\QQ/\ZZ$, and hence the quotient must
vanish. Consequently, $D+{\operatorname{im}}\rmapt_A=\pKz{A}$.

As $\ker\bmapt_A={\operatorname{im}}\rmapt_A$, we infer that $\bmapt_A\restr{D}$ is an
isomorphism, and we let $\sigma=(\bmapt_A\restr{D})\inv$. As
$\bmapt_A\tau=\bmapt_I\tau=\id$ and ${\operatorname{im}}\tau\subseteq D$, $\sigma$ does
extend $\tau$. 
\end{pf}

Combining 6.1 and 8.1 of \cite{gae:dgt}, we get that the lattice
isomorphism between ideals of $A$ and order ideals of $K_0(A)$
extends to a lattice isomorphism
\[
I\mapsto K_0(I)\oplus K_1(I)
\]
into the order ideals of $K_0(A)\oplus K_1(A)$ for the
\cstar-algebras we consider. As this is also a Riesz group, sums and
intersections of order ideals are again order ideals. We conclude, in
particular, that
\[
K_1(I\cap J)=K_1(I)\cap K_1(J)\qquad
K_1(I + J)=K_1(I)+ K_1(J)
\]

Let $G_1,\dots,G_n$ be a set of subgroups of a group $H$. There are
maps
\[
\Gamma^1:\bigoplus_{i<j}{G_i\cap G_j}\rightarrow \bigoplus_i{G_i}
\qquad
\Gamma^0: \bigoplus_i{G_i}\rightarrow H
\]
given on each summand by $\Gamma^0(g_i)=g_i$ and
\[
\Gamma^1(g_{ij})_k=\left\{\begin{array}{ll}
                          g_{ij}&k=i\\-g_{ij}&k=j\\0&\text{other }k
                          \end{array}\right.
\]
Note that $\Gamma^0\Gamma^1=0$.

\begin{lemma}\label{maxideals}
Suppose $I$ and $J$ are ideals in a real rank
zero $AD$ algebra $A$ such that $I+J=A$. Then
$$
0\longrightarrow
{K_1(I\cap J)}\stackrel{\Gamma^1}{\longrightarrow}
{K_1(I)\oplus K_1(J)}\stackrel{\Gamma^0}{\longrightarrow}
{K_1(A)}\longrightarrow 0
$$
is a pure exact extension.
\end{lemma}
\begin{pf}
Assume that $n(x_1,x_2)$ lies in the image of $\Gamma^1$. There is
then $x\in K_1(I\cap J)$ with $x=nx_1$. As $I\cap
J$ is an ideal of $I$, the subgroup  $K_1(I\cap J)$ is pure in
$K_1(I)$, which implies that $x=ny$ with $y\in K_1(I\cap J)$. Clearly
$n(x_1,x_2)=n\Gamma^1(y)$.
\end{pf}

The following result is very similar to \cite[4.11]{gg:ccrrzuet} and is
proved by the exact same argument. We include it for the convenience
of the reader.

\begin{lemma}\label{exact}
When $(I_j)$ is a finite family of comaximal ideals of a real rank
zero $AD$ algebra, the complex
$$
{\bigoplus_{i<j}{\tor K_1(I_i\cap I_j)}}\stackrel{\Gamma^1}{\longrightarrow}
{\bigoplus_i{\tor K_1(I_i)}}\stackrel{\Gamma^0}{\longrightarrow}
{\tor K_1(A)}\longrightarrow 0
$$
is exact.
\end{lemma}
\begin{pf}
Surjectivity of $\Gamma^0$ follows from Lemma \ref{maxideals}.
We prove exactness at $\bigoplus_i{\tor K_1(I_i)}$ for $n=2$ and
$n=3$. A straightforward induction argument will then prove the
general claim. The case $n=2$ is obvious from Lemma \ref{maxideals}.
For the case $n=3$, assume that
$$
y_1+y_2+y_3=0
$$
with $y_i\in\tor K_1(I_i)$. As the lattice of ideals in $A$ is
distributive (since $I^2=I$ for every ideal $I$), we have
$$
I_1=I_1\cap(I_2+I_3)=(I_1\cap I_2)+(I_1\cap I_3),
$$
and we apply Lemma \ref{maxideals} to $I_1$ to get that
$y_1=z_2+z_3$
for some $z_i\in\tor K_1(I_1\cap I_i)$. From the claim for $n=2$ we get
that
$$
(0,y_2+z_2,y_3+z_3)=(0,w,-w)
$$
for some $w\in\tor K_1(I_2\cap I_3)$. And then
$$
(y_1,y_2,y_3)=\Gamma^1(z_2,z_3,w).
$$
\end{pf}

{\bf Proof of 4.1:}\label{existidsplit}
Let $\Omega$ be a proper, hereditary subset of the lattice of ideals
in $A$. We get from the finiteness assumption that there is an ideal
$I\not\in\Omega$ with $\min(I)$, the set of proper ideals of $I$,
contained in $\Omega$. This means that $\Omega\cup\{I\}$ is again
hereditary, and we can prove the claim inductively by showing that
when $\Omega$ is such a set, and K\"unneth splittings $\sigma_J$ are
given for all $J\in\Omega$ with
\begin{eqnarray}\label{cohere}
J_1\subseteq J_2, J_i\in \Omega\infer \sigma_{J_2}\restr{\tor
K_1(J_1)}=\sigma_{J_1},
\end{eqnarray}
and $I\not\in\Omega$ is given with $\min(I)\subseteq \Omega$, then we
can define $\sigma_I$ such that \eqref{cohere} holds true for
$\Omega\cup\{I\}$. Let $N$ be the number of maximal elements of
$\min(I)$. We must deal with the cases $N=1$ and $N>1$ separately.

\gritem{1} $N=1$: Denote the single maximal ideal by $J$ ($J$ could be
the zero ideal). By Lemma \ref{maximalcase}, $\sigma_J$ extends to a
splitting map $\sigma_I$, and \eqref{cohere} follows from
\[
\sigma_{I}\restr{\tor K_1(J_1)}=
\sigma_{J}\restr{\tor K_1(J_1)}=\sigma_{J_1}.
\]

\gritem{2} $N>1$: Let $I_1,\dots I_N$ be the set of different maximal
ideals. Clearly $I_i+I_j=I$ for all $i\not=j$, and we get a diagram
\begin{flushleft}
\begin{picture}(00,140)
\put(100,130){$0$}
\put(220,130){$0$}
\put(100,105){$\vector(-0,1){20}$}
\put(220,105){$\vector(-0,1){20}$}
\put(80,90){${\pKz{I}}$}
\put(210,90){${\tor K_1(I)}$}
\put(200,90){$\vector(-1,0){50}$}
\put(70,50){${\displaystyle\bigoplus_i{\pKz{I_i}}}$}
\put(100,60){$\vector(-0,1){20}$}
\put(220,60){$\vector(-0,1){20}$}
\put(230,70){${\Gamma}_0$}
\put(85,70){${\Gamma}_0$}
\put(200,50){${\displaystyle\bigoplus_i{\tor K_1(I_i)}}$}
\put(190,55){$\vector(-1,0){40}$}
\put(165,58){${\oplus{\sigma_{i}}}$}
\put(60,10){${\displaystyle\bigoplus_{i<j}{\pKz{I_i\cap I_j}}}$}
\put(100,20){$\vector(-0,1){20}$}
\put(220,20){$\vector(-0,1){20}$}
\put(85,30){${\Gamma}_1$}
\put(200,10){${\displaystyle\bigoplus_{i<j}{\tor K_1(I_i\cap I_j)}}$}
\put(230,30){${\Gamma}_1$}
\put(190,10){$\vector(-1,0){30}$}
\put(165,15){${\oplus{\sigma_{ij}}}$}
\end{picture}
\end{flushleft}
\noindent
in which the vertical maps form a complex and the solid square is
commutative by \eqref{cohere}. By Lemma \ref{exact},
the horizontal maps induce a map $\sigma_I:\tor K_1(I)\rightarrow \pKz
{I}$, such that the other square commutes. This implies that
\eqref{cohere} holds for $\Omega\cup\{I\}$ also.
\begin{flushright}$\Box$\end{flushright}
\begin{remark}\label{finalwords}
The method of proof given here differs from the one given in
\cite{gg:ccrrzuet} by us being able to argue by sheer algebraic
means, and without 
going back into the inductive limit presentation of the \cstar-algebra
in question. That is of course only possible because we have a
complete algebraic invariant at our service.
On the other hand, the reader will note that the overall structure of
the two proofs are the same, in particular in the way we deal with the
lattice of ideals. We have certainly been inspired by (cite{gg:ccrrzuet}.
\end{remark}

\begin{remark}\label{ASH}
When working with general $ASH$ algebras of real rank zero and with
slow dimension growth, the invariant consists of
the full collection of groups $\pKsubi{n}{A}$ along with order
structures and coherence maps as described in \linebreak \DG{4}. 
As mentioned above, every unspliced sequence
\[
0\longrightarrow{K_i(A)\otimes\ZZp{n}}\stackrel{\rmapt_n^i}{\longrightarrow}{\pKsubi{n}{A}}\stackrel{\bmapt_n^i}{\longrightarrow}
{K_{i+1}(A)[n]}\longrightarrow 0
\]
is split. In fact, as in
\cite{cfb:skskI},\cite{cfb:skskII} one can choose an entire family
of splitting maps which is {\em coherent}\/ in the coefficient.  This
means, when we denote by $\lmapk{m}{n}^i$ the natural
map
\[ 
\lmapk{m}{n}^i:K_i(A)[n]\rightarrow K_i(A)[m], \] that
$\sigma_m^i\lmapk{m}{n}^{i+1}=\kmapk{m}{n}^i\sigma_n^i$. Such a
coherent family is what we understand by a K\"unneth splitting in this
setting.

We define {\it ideal} K\"unneth splittings as above, and get the natural analogue of
Proposition \ref{idealsplit} from \DG{9.1-2}. We can prove that $ASH$
algebras of real rank zero with
finitely many ideals are ideally split by proceeding as in the proof
of Proposition  \ref{existidspli}. Here, our approach in case $2^\circ$ easily generalizes
because
\[
{\bigoplus_{j<k}{K_i(I_j\cap I_k)[n]}}\stackrel{\Gamma^1}{\longrightarrow}
{\bigoplus_j{K_i(I_j)[n]}}\stackrel{\Gamma^0}{\longrightarrow}
{K_i(A)[n]}\longrightarrow 0
\]
is exact as a consequence of purity, but we only know of a fairly
laborious way to generalize our approach in case $1^\circ$. A possible
method is to achieve an analogue of
Lemma \ref{maximalcase}, by hands-on extending a given coherent 
family defined on the $K$-theory of an ideal. This is seen by
inspection of the proofs of 
\cite[2.8]{cfb:skskI}, \cite[2]{cfb:skskII}, using at a crucial point
that if a basis (cf. \cite[17.2]{lf:iagI}) for $K_1(I)[p^r]$ is given, it
can be augmented to a basis for  $K_1(A)[p^r]$ by purity arguments.
\end{remark}

\begin{center}
{\bf Acknowledgements}
\end{center}

I am grateful to Peter Friis and Niels Peter J{\o}rgensen for helpful discussions, and to Ken Goodearl for comments on the first version of the paper which lead to a substantial increase in clarity and a dramatic reduction of the proof of Lemma 4.2.

\providecommand{\bysame}{\leavevmode\hbox to3em{\hrulefill}\thinspace}
\providecommand{\MR}{\relax\ifhmode\unskip\space\fi MR }
\providecommand{\MRhref}[2]{%
  \href{http://www.ams.org/mathscinet-getitem?mr=#1}{#2}
}
\providecommand{\href}[2]{#2}

\end{document}